\documentclass[12pt]{article}
\usepackage{makeidx}
\usepackage{xcolor}
\usepackage{amsfonts,amsmath,amsthm}
\usepackage{wrapfig}
\usepackage{amstext}
\usepackage{amssymb}
\usepackage{epsfig}
\usepackage{color}
\usepackage{graphicx}
\usepackage{graphics}
\usepackage{hyperref}

\def \squarebox#1{\hbox to #1{\hfill \vbox to #1{\vfill}}}
\newtheorem{teorema}{Theorem}[section]
\newtheorem{lema}[teorema]{Lemma}
\newtheorem{corolario}[teorema]{Corollary}
\newtheorem{proposicao}[teorema]{Proposition}

\newtheorem*{teoremaNo}{Theorem}

\newtheorem{teodefi}[teorema]{Theorem -- Definition}
\newenvironment{prova}{\noindent {\bf Proof:}}{\hfill $\qed $ \newline}
\newenvironment{exemplo}{\vspace{.3cm}\noindent {\bf Example:}}{\hfill $\qed$ \newline}
\newtheorem{remark}[teorema]{Remark}

\DeclareMathOperator{\ad}{ad}  
\DeclareMathOperator{\Ad}{Ad}  
  
\DeclareMathOperator{\tr}{tr}

\begin{document}

\title{Semigroups in semi-simple Lie groups: Flag type and estimation of cocycles}
\author{Luiz A.B. San Martin\thanks{Supported by CNPq grant no.\ 303755/09-1, FAPESP grant no.\
2012/18780-0 and CNPq/Universal grant no 476024/2012-9.} and Victor Uzita \\
		Instituto de Matem\'{a}tica\\
		Universidade Estadual de Campinas, Brazil \and
Adriano Da Silva\thanks{Supported by Proyecto UTA Mayor Nº 4781-24}\\
		Departamento de Matem\'atica,\\Universidad de Tarapac\'a - Iquique, Chile.}
	\date{\today }
	\maketitle
	
	\begin{abstract}
		The flag type of a semigroup $S$ of a noncompact semisimple Lie group is an algebraic tool related to the geometry of the invariant control set determined by $S$ on the flag manifolds of $G$. In the present paper we show that it is possible to recover the flag type by studying the existence of lower bounds for cocycles on the maximal flag manifold. 
	\end{abstract}
	
	

\section{Introduction}

Let $G$ be a noncompact semisimple Lie group $G$ with a finite center and $
S\subset G$ a subsemigroup with nonempty interior. The purpose of this paper
is to characterize the \textit{flag type} of $S$ by means of lower bounds of
cocycles over the flag manifolds of $G$.

The concept of the flag type of a semigroup has proved to be valuable for the
study of semigroups in semisimple Lie groups (see \cite{smconex, smorder, smt, smsan}) and its applications to dynamical
systems and harmonic analysis (see \cite{smpoiss, patrsm, smcone, smlucas, smmoment}). It appears in the study of the action of a semigroup $S$ on the flag manifolds of $G$. Precisely, we have the following result, proved in \cite[Theorem 4.3]
{smt}, ensuring its existence and uniqueness (see also \cite{smmax, smc, smorder}).
\begin{teodefi}
\label{flagtype}
Let $S\subset G$ be a proper semigroup with $\mathrm{int}S\neq \emptyset $.
Then there exists a unique flag manifold $\mathbb{F}_{\Theta \left( S\right)
}$ satisfying the following two conditions:

\begin{enumerate}
\item $\mathbb{F}_{\Theta \left( S\right) }$ is maximal among the flag
manifolds $\mathbb{F}_{\Theta }$ such that the unique $S$-invariant control
set $C_{\Theta }\subset \mathbb{F}_{\Theta }$ is contractible in the sense
that there exists $g\in \mathrm{int}S$ such that  $g^{n}C_{\Theta}$ shrinks to a point as $n\rightarrow +\infty$.

\item $\mathbb{F}_{\Theta \left( S\right) }$ is minimal among the flag
manifolds $\mathbb{F}_{\Theta }$ such that $C=\pi ^{-1}\left( C_{\Theta
}\right) $ is the invariant control set in the maximal flag manifold $%
\mathbb{F}$.
\end{enumerate}

The flag manifold $\mathbb{F}_{\Theta \left( S\right) }$ and the subset of roots $\Theta(S)$ are called, indistinguishly, the 
\emph{flag type} of the semigroup $S$.
\end{teodefi}

All the conditions for the flag type of $S$ use in an essential way that $S$ has
nonempty interior.

In this paper we prove another equivalent condition that relates the flag
type to lower bounds of cocycles. The statement of this condition does not
require in advance that the semigroup has a nonempty interior and hence makes
sense to other classes of semigroups.

Let $G=KAN$ be an Iwasawa decomposition. By this decomposition the minimal parabolic
subgroup is $P=MAN$ where $M$ is the centralizer of $A$ in $K$ and $\mathbb{F}=G/P=K/M$. The $K$-invariant cocycles over the flag manifolds are defined after this decomposition as follows: Define the map $\rho :G\times K\rightarrow A$ by%
\begin{equation*}
gu=k\rho \left( g,u\right) n\in KAN
\end{equation*}%
and put $\mathsf{a}\left( g,u\right) =\log \rho \left( g,u\right) \in
\mathfrak{a}$. These maps are right $M$-invariant in the second variable
hence they factor to maps (with the same notation) $\rho :G\times \mathbb{F}\rightarrow A$ and $\mathsf{a}:G\times \mathbb{F}%
\rightarrow \mathfrak{a}$. For $%
\lambda \in \mathfrak{a}^{\ast }$ we write $\mathsf{a}_{\lambda }\left(
g,x\right) =\lambda \mathsf{a}\left( g,x\right) $ and $\rho _{\lambda
}\left( g,x\right) =e^{\mathsf{a}_{\lambda }\left( g,x\right) }$.

Next we state the main result of this paper (Theorem \ref{teoprinc}).

\begin{teoremaNo}
Let $S\subset G$ be a semigroup with $\mathrm{int}S\neq
\emptyset $ and write $\mathbb{F}_{\Theta \left( S\right) }$, $\Theta \left(
S\right) \subset \Sigma $, for its flag type. Denote by $C$ the invariant
control set of $S$ in the maximal flag manifold $\mathbb{F}$.

If $x_{0}\in C_{0}$, then 
\begin{equation*}
\inf_{g\in S}\rho _{\alpha }\left( g,x_{0}\right) >0
\end{equation*}%
for any $\alpha \in \Sigma \setminus \Theta \left( S\right) $. Conversely, $%
\inf_{g\in S}\rho _{\alpha }\left( g,x_{0}\right) =0$ if $\alpha \in \Theta
\left( S\right) $.
\end{teoremaNo}

Let us mention that, the existence of a positive lower bound for cocycles on flag manifolds has already been used in \cite{smmoment} in the study of the moment Lyapunov exponents of the i.i.d. random product on $G$ defined by a probability measure. This leads to several important results concerning the spectral radii of compact operators on Banach spaces. Another importance of the previous result, is that for the statement of the conditions in this theorem, it is
not required that the semigroup has a nonempty interior. The condition can be
stated for any semigroup that has a unique invariant control set in $\mathbb{
F}$. Thus, apart from its intrinsic interest for the structure of the semigroups
with nonempty interior, the above theorem can open the way to the understanding of semigroups in a more general context. We have in mind the completely irreducible and contracting semigroups that
arise in the study of random products in semi-simple Lie groups (see {Guivarc'h-Raugi }\cite{gr}, Abels-Margulis-Soifer \cite{abelmar}, and
Gol'dsheid-Margulis \cite{goldmarg}). In case $G$ is an algebraic group, irreducibility plus contractibility is equivalent to saying that the semigroup
is Zariski dense.

The paper is structured as follows: In Section 2 we obtain lower bounds for the cocycles where the associated functional belongs to the cone generated by fundamental weights outside the flag type. In particular, a very useful result for rank-one groups is obtained. Such a result allows us to restrict the proof of our main result to the fibers of some specific fibration between flag manifolds. We finish the section analyzing when such a lower bound is uniform. By a concrete example, we show that uniformity does not hold in general; however, one can find a coset of the semigroup where uniformity holds.

In Section 3 we state and prove our main result relating the flag type with the lower bounds of cocycles on the flag manifolds. We start by showing that the cocycles associated with root in the flag type of the semigroup admit no lower bound. In sequence, we look at roots outside the flag type. By looking at a fibration of the maximal flag onto the flag associated with the root, we are able to show that the values of the cocycle, for a fixed point in the core of the invariant control set, coincide with a cocycle defined on the invariant control set of a rank-one group. This fact, together with the results in Section 3, allows us to obtain the desired lower bound. In order to make the paper self-contained and more fluent, an appendix is available, where we introduced the basic results and notations related to semisimple theory. We also use such a section to define more formally the flag type of semigroups and the cocycles induced by the Iwasawa decomposition on the flag manifolds.

\section{Lower bounds of cocycles}

\label{seclower}

In this section we present a series of lemmas which will help us to prove our main results. Although some of these lemmas were already proved in \cite{smmoment}, in order to keep the paper self-contained we reproduce their proof here again

Let us, as previously, consider $S\subset G$ be a semigroup wit $\mathrm{int}S\neq \emptyset $ and denote
by $C$ the unique invariant control set in the maximal flag manifold. In a
partial flag manifold $\mathbb{F}_{\Theta }$ the invariant control set is
denoted by $C_{\Theta }$.

For $x$ in a control set we write $S_{x}=\{g\in S:gx=x\}$ and $S_{x}^{\circ
}=\{g\in \mathrm{int}S:gx=x\}$. The core of a control set is the subset of
those $x$ such that $S_{x}^{\circ }\neq \emptyset $. The core of the control set $C_{\Theta }$ is denoted by $(C_{\Theta}) _{0}$.

The proof of the lower bound in Theorem \ref{teoprinc} is based in the
following lemma that reduces the estimate to the isotropy group.

\begin{lema}
\label{lemcompa}Let $S$ be a semigroup with $\mathrm{int}S\neq \emptyset $
and denote by $C$ its invariant control set in the maximal flag manifold $%
\mathbb{F}=G/P$. Take $x_0\in C_{0}$ and $\lambda \in \mathfrak{a}^{\ast }$
and suppose that there exists $d>0$ such that $\rho _{\lambda }\left(
g,x_0\right) >d$ for all $g\in S^{\circ}_{x_0}$.

Then there exists $c>0$ such that $\rho _{\lambda }\left( g,x_0\right) >c$ for
all $g\in S$.
\end{lema}

\begin{prova}
Suppose by contradiction the existence of a sequence $g_{k}\in S$ with $%
\rho _{\lambda }\left( g_{k},x_0\right) \rightarrow 0$. Since $C$ is a compact and invariant subset, it can be assumed that 
$g_{k}x_0\rightarrow y$ with $y\in C$. If $g\in S$ then 
$$\rho _{\lambda }(gg_{k},x_0)=\rho _{\lambda
}(g,g_{k}\cdot x_0)\rho _{\lambda }(g_{k},x_0)\rightarrow 0,$$
because the
function $z\mapsto \rho _{\lambda }\left( g,z\right) $ is bounded. Take in
particular $g\in \mathrm{int}S$ with $gy=x_0$. Then we can substitute $g_{k}$
by $gg_{k}$ and assume that $g_{k}\in \mathrm{int}S$ and $g_{k}x\rightarrow
x_0 $.

Since $x_0\in C_{0}$ the semigroup $S_{x_0}^{\circ }$ is not empty. Take $%
g_{0}\in S_{x_0}^{\circ }$ and a compact neighborhood $W$ of $g_{0}$ in $%
\mathrm{int}S$. We have that $U=W^{-1}x_0$ is a neighborhood of $x_0$ in $%
\mathbb{F}$ because $g_{0}x_0=x_0$. By construction, for every $z\in U$ there
exists $h\in W$ such that $x_0=hz$. Write 
\begin{equation*}
r=\sup \{\rho _{\lambda }\left( h,z\right) :h\in W,z\in \mathbb{F}\}
\end{equation*}%
which is finite by compactness.

Now, let $k$ be large enough so that $g_{k}x_0\in U$ and $\rho _{\lambda
}\left( g_{k},x_0\right) <d/2r$. Then there exists $h\in W$ such that $%
hg_{k}x_0=x_0$ and we have 
\begin{equation*}
\rho _{\lambda }\left( hg_{k},x_0\right) =\rho _{\lambda }(h,g_{k}x_0)\rho
_{\lambda }(g_{k},x_0)\leq r\rho _{\lambda }\left( g_{k},x_0\right) <\frac{rd}{2r%
}=\frac{d}{2}
\end{equation*}%
contradicting the assumption.
\end{prova}

To get a first application of Lemma \ref{lemcompa}, we introduce the
following notation: Write the simple system of roots as $\Sigma =\{\alpha
_{1},\ldots ,\alpha _{l}\}$ and let $\Phi =\{\mu _{1},\ldots ,\mu
_{l}\}\subset \mathfrak{a}^*$ be the set of corresponding fundamental weights, that is, the linear functionals satisfying
\begin{equation*}
\frac{2\langle \alpha _{i},\mu _{j}\rangle }{\langle \alpha _{i},\alpha
_{i}\rangle }=\delta _{ij}.
\end{equation*}%
Take $\Theta \subset \Sigma $ with $\Theta =\{\alpha _{i_{1}},\ldots ,\alpha
_{i_{j}}\}$ and let $\Phi _{\Theta }=\{\mu _{i_{1}},\ldots ,\mu
_{i_{j}}\}$ be the set of fundamental weights with the same indices as those
in $\Theta $. Equivalently 
\begin{equation*}
\Phi \setminus \Phi _{\Theta }=\{\mu \in \Phi :\forall \alpha \in \Theta
,~\langle \alpha ,\mu \rangle =0\}.
\end{equation*}%
Denote by $\left( \mathfrak{a}_{\Theta }^{\ast }\right) ^{+}$ the
\textquotedblleft partial chamber\textquotedblright\ 
\begin{equation*}
\left( \mathfrak{a}_{\Theta }^{\ast }\right) ^{+}=\{\beta \in \mathfrak{a}%
_{\Theta }^{\ast }:\forall \alpha \in \Sigma \setminus \Theta ,~\langle
\alpha ,\beta \rangle >0\}
\end{equation*}%
which is the interior (in $\mathfrak{a}_{\Theta }^{\ast }$) of the convex
cone $\mathrm{cl}\left( \mathfrak{a}_{\Theta }^{\ast }\right) ^{+}$ spanned
by $\Phi \setminus \Phi _{\Theta }$.

Now assume without loss of generality that $A^{+}\cap \mathrm{int}S\neq
\emptyset $ where $A^{+}=\exp \mathfrak{a}^{+}$. This assumption implies
that origin $x_{0}$ of $\mathbb{F}=G/MAN$ belongs to the core $C_{0}$ of the
invariant control set $C$. The convex cone 
\begin{equation*}
\Gamma _{N}=\{H\in \mathfrak{a}:\exists n\in N,\exists t>0,~e^{tH}n\in 
\mathrm{int}S\}
\end{equation*}%
was considered in \cite[Section 4]{smt}. The following inclusion was
proved there.

\begin{proposicao}
\label{propConesEmA}Let $\mathbb{F}_{\Theta }$, $\Theta \subset \Sigma $ be
the flag type of $S$. Then 
\begin{equation}
\Gamma _{N}\subset \Gamma _{\Theta }=\mathrm{cl}\left( \bigcup\limits_{w\in 
\mathcal{W}_{\Theta }}w\mathfrak{a}^{+}\right) .  \label{forIncluConeA}
\end{equation}%
Moreover, $\Gamma _{N}\cap w\mathfrak{a}^{+}\neq \emptyset $ for every $w\in 
\mathcal{W}_{\Theta }$.
\end{proposicao}

\begin{remark}
    In \cite{smt} it is proved that if the inclusion fails, then there are $%
w\notin \mathcal{W}_{\Theta }$, $H\in w\mathfrak{a}^{+}$ and $n\in N$ such
that $g=e^{H}n\in \mathrm{int}S$. But this implies that $x_{0}$ is a fixed
point of type $w$ of $g$. Since $w\notin \mathcal{W}_{\Theta }$ this
contradicts the fact that $\Theta $ is the flag type of $S$.
\end{remark}

From the inclusion (\ref{forIncluConeA}) we can prove the following estimate
for elements in $\mathrm{int}S$ fixing a point $x_{0}\in \mathbb{F}$.

\begin{lema}
\label{lemacocmaior1}Suppose that $g\in \mathrm{int}S$ is such that $%
gx_{0}=x_{0}$ and take $\lambda \in \left( \mathfrak{a}_{\Theta }^{\ast
}\right)^{+}$, where $\Theta$ is the flag type of $S$. Then $\rho _{\lambda }(g,x_{0})\geq 1$.
\end{lema}

\begin{prova}
If $gx_{0}=x_{0}$, then we can write $g=man\in MAN$ in which case $\rho
_{\lambda }\left( g,x_{0}\right) =e^{\lambda \left( H\right) }$ where $%
a=e^{H}$. We claim that $H\in \Gamma _{\Theta }$. In fact, since the elements of finite order are dense in compact groups, we can perturb $g$
inside $\mathrm{int}S$ and assume that $m$ has finite order. Then for some $%
j\geq 1$ we have $g^{j}=a^{j}\overline{n}\in AN$, which implies, by the
comments above, that $jH\in \Gamma _{N}\subset \Gamma _{\Theta }$ and hence $%
H\in \Gamma _{\Theta }$. Finally, for any $\mu\in\Phi\setminus\Phi_{\Theta}$, $\alpha\in\Theta$ and $H\in\mathfrak{a}$, it holds that $\mu(r_{\alpha}(H))=\mu(H)$. Since any $\left(\mathfrak{a}%
_{\Theta }^{\ast }\right) ^{+}$ is generated by $\Phi \setminus \Phi _{\Theta }$, the previous implies, in particular, that $\lambda(H)\geq 0$ and hence $\rho
_{\lambda }\left( g,x_{0}\right) =e^{\lambda \left( H\right) }\geq 1$ as
claimed.
\end{prova}

Combining this lemma with Lemma \ref{lemcompa}, we get at once the following
estimate of cocycles.

\begin{proposicao}
\label{propLowerCone}Let $S$ be a semigroup whose flag type is $\mathbb{F}%
_{\Theta }$, $\Theta \subset \Sigma $. Denote by $C$ its invariant control
set in the maximal flag manifold $\mathbb{F}$. Take $\lambda \in \left(
\mathfrak{a}_{\Theta }^{\ast }\right) ^{+}$ and $x\in C_{0}$. Then there
exists $c>0$ such that $\rho _{\lambda }\left( g,x\right) >c$ for all $g\in
S $.
\end{proposicao}

This proposition will be applied later to get the lower bound estimate of
Theorem \ref{teoprinc}. For this application it is required only the
following specialization to real rank one groups.

\begin{corolario}
\label{corLowerRankone}Suppose that $\mathfrak{g}$ (and $G$) has real rank
one and let $S\subset G$ be a proper semigroup with $\mathrm{int}S\neq
\emptyset $. Denote by $\alpha $ (and eventually $2\alpha $) the positive
root. Then $\inf_{g\in S}\rho _{\alpha }\left( g,x\right) >0$ if $x\in C_{0}$%
.
\end{corolario}

\begin{prova}
In the rank one case there is only the maximal flag manifold $\mathbb{F}$
and if $S$ is proper, then its flag type is $\mathbb{F}=\mathbb{F}_{\emptyset
}$ itself. The subspace $\mathfrak{a}_{\Theta }^{\ast }$ is one-dimensional.
and $\left( \mathfrak{a}_{\Theta }^{\ast }\right) ^{+}$ is the ray
containing $\alpha $. Hence $\alpha $ falls in the condition of the above
proposition so that $\rho _{\alpha }\left( g,x\right) $ is lower bounded.
\end{prova}

The following example illustrates this corollary with a concrete semigroup
in $\mathrm{Sl}\left( 2,\mathbb{R}\right) $.
%TCIMACRO{\TeXButton{vspace12pt}{\vspace{12pt}}}%
%BeginExpansion
\vspace{12pt}%
%EndExpansion

%TCIMACRO{\TeXButton{noindent}{\noindent}}%
%BeginExpansion
\noindent%
%EndExpansion
\begin{exemplo}
For $\mathrm{Sl}\left( 2,\mathbb{R}\right) $ the only flag
manifold is the projective line $\mathbb{P}^{1}$. The cocycle $\rho
_{\lambda }$ over $\mathbb{P}^{1}$ is defined by the relation 
\begin{equation*}
\lambda \left( 
\begin{array}{cc}
1 & 0 \\ 
0 & -1%
\end{array}%
\right) =1
\end{equation*}%
and hence, $\rho _{\lambda }\left( g,\left[ z\right] \right) =\left\Vert
gz\right\Vert /\left\Vert z\right\Vert $, $0\neq z\in \mathbb{R}^{2}$.
Consider the cone 
\begin{equation*}
W=\{\left( a,b\right) \in \mathbb{R}^{2}:a\geq 0,~|b|\leq a\}
\end{equation*}%
and define the semigroup $S_{W}=\{g\in \mathrm{Sl}\left( 2,\mathbb{%
R}\right) :gW\subset W\}$. The core of the invariant control set for the action of $S_W$ in $\mathbb{%
P}^{1}$ is%
\begin{equation*}
C_{0}=\{\left[ \left( a,b\right) \right] \in \mathbb{P}^{1}:\left(
a,b\right) \in \mathrm{int}W\}.
\end{equation*}

Let $g\in S_{W}$ satisfying $g\left( 1,0\right) =\left( a,b\right) \in W$. Then, for 
\begin{equation*}
h=\left( 
\begin{array}{cc}
1 & 0 \\ 
b/a & 1%
\end{array}%
\right) \left( 
\begin{array}{cc}
a & 0 \\ 
0 & a^{-1}%
\end{array}%
\right),
\end{equation*}%
it holds that that $g\left( 1,0\right) =h\left( 1,0\right)$. Therefore,  $h^{-1}g\left( 1,0\right) =\left( 1,0\right) $ implying that $h^{-1}g$ is upper
triangular, and hence, $g$ has the form%
\begin{equation*}
g=\left( 
\begin{array}{cc}
1 & 0 \\ 
y & 1%
\end{array}%
\right) \left( 
\begin{array}{cc}
\mu & 0 \\ 
0 & \mu ^{-1}%
\end{array}%
\right) \left( 
\begin{array}{cc}
1 & x \\ 
0 & 1%
\end{array}%
\right) .
\end{equation*}%
with $\mu >0$ and $|y|\leq 1$. Since, $g\left( 1,-1\right) \in W$
and $g\left( 1,1\right) \in W$, it holds that  
$$1-x\geq 0, \hspace{.5cm}|\mu
y\left( 1-x\right) -\mu ^{-1}|\leq \mu \left( 1-x\right), $$
$$1+x\geq 0 \hspace{.5cm} \mbox{ and }\hspace{.5cm}
|\mu y\left( 1-x\right) +\mu ^{-1}|\leq \mu \left( 1+x\right). $$ 
Hence 
$$
|x|\leq 1  \hspace{.5cm} \mbox{ and }\hspace{.5cm}
\mu ^{-1}-\mu |y|\left( 1-x\right) \leq \mu \left( 1-x\right),$$
that is, $\mu ^{-1}\leq \mu \left( 1-x\right) \left( 1+|y|\right) $. But $%
|y|\leq 1$ and $1-x\leq 2$ so that $1\leq 4\mu ^{2}$ and $\mu \geq 1/2$.
Since $g\left( 1,0\right) =\mu \left( 1,y\right) $ we get the lower bound 
\begin{equation*}
\left\Vert g\left( 1,0\right) \right\Vert =\mu \sqrt{1+y^{2}}\geq \frac{1}{2}.
\end{equation*}
Now if $z=\left[ \left( a,b\right) \right] \in C_{0}$, there exists $%
h\in S_{W}$ satisfying $z=hz_{0}$, where  $z_{0}=\left[ \left( 1,0\right) \right] $.
Consequently, if $g\in S_{W}$ then $gh\in S_{W}$ so that 
\begin{equation*}
\rho _{\lambda }\left( g,z\right) =\frac{\left\Vert ghz_{0}\right\Vert }{%
\left\Vert hz_{0}\right\Vert }\geq \frac{\left\Vert ghz_{0}\right\Vert }{%
\left\Vert h\right\Vert }\geq \frac{1}{2\left\Vert h\right\Vert }, 
\end{equation*}
showing that $c=1/2\left\Vert h\right\Vert $ is the desired lower bound for $z=hz_{0}\in C_{0}$.
Notice that the lower bound $c=1/2\left\Vert h\right\Vert $
depends on $z$ through $\left\Vert h\right\Vert $ so that it is not uniform
in $z$. Actually, a uniform lower bound (in $C_{0}$) cannot be obtained. In
fact
\begin{equation*}
h_{t}=\left(
\begin{array}{cc}
\cosh t & \sinh t \\
\sinh t & \cosh t%
\end{array}%
\right)
\end{equation*}%
belongs to $S_{W}$ for all $t\in \mathbb{R}$. If $z=\left[ \left( a,b\right) %
\right] \in C_{0}$ then
\begin{equation*}
\rho _{\lambda }\left( h_{t},z\right) =\frac{\left( a^{2}+b^{2}\right) \cosh
2t+2ab\sinh 2t}{a^{2}+b^{2}}
\end{equation*}%
which for fixed $t$ converges to $e^{-2t}$ as $\left( a,b\right) $ approaches
$\left( 1,-1\right) $. This shows that a uniform lower bound does not necessarily exist.

\end{exemplo}

The previous example shows that the lower bound $c$ in Proposition \ref{propLowerCone} could depend on
$x$ and may not be uniform. Despite this fact, the next proposition shows that a uniform lower bound can be found for some coset $Sh$.

\begin{proposicao}
Under the assumptions of Proposition \ref{propLowerCone}, there are $h\in S$ and $c>0$
such that $\rho _{\lambda }\left( g,y\right) >c$ for all $g\in Sh$ and $y\in
C$.
\end{proposicao}

\begin{prova} Let $x_0\in C_0$ and consider $g_0\in\mathrm{int}S$ such that $g_0x_0=x_0$. Let $U\subset\mathrm{int}S$ a compact neighborhood of $g_0$ and write $W=U\xi_0$, where $\xi_0:=\pi(x_0)$, with $\pi:\mathbb{F}\rightarrow\mathbb{F}_{\Theta}$ the canonical projection. Note that 
$W$ is a compact neighborhood of $\xi_0$ in $\mathbb{F}_{\Theta}$, since $$\xi_0=\pi(x_0)=\pi(g_0x_0)=g_0\pi(x_0)=g_0\xi_0\in W.$$
Moreover, 
$$x_0\in C_0\hspace{.5cm}\implies\hspace{.5cm}\xi_0\in(C_{\Theta})_0\hspace{.5cm}\implies\hspace{.5cm} W\subset (C_{\Theta})_0,$$
where the last inclusion follows by the $S$-invariance of $C_{\Theta}$. By Theorem \ref{flagtype}, there exists $n_0\in\mathbb{N}$ such that $h:=g_0^{n_0}$ satisfies $hC_{\Theta}\subset W$ and hence $\pi^{-1}(W)$ is a compact neighborhood of $x_0$ satisfying, 
$$hC=h\pi^{-1}(C_{\Theta})=\pi^{-1}(hC_{\Theta})\subset \pi^{-1}(W)\subset\pi^{-1}(C_{\Theta})=C.$$

Since $\lambda\in (\mathfrak{a}_{\Theta}^*)^+$, we get by Lemma \ref{restriction} and Proposition \ref{propLowerCone} that
$$\exists c_1>0, \hspace{.5cm}\rho_{\lambda}(g, ux_0)=\rho_{\lambda}(g, x_0)>c_1, \hspace{.5cm}\mbox{ for all } \hspace{.5cm} g\in S,  u\in K_{\Theta}.$$
Therefore, for any $g'\in U$ and $u\in K_{\Theta}$ we get
\begin{equation*}
\rho _{\lambda }\left(gg', ux_0\right) =\rho _{\lambda }\left(g, g'ux_0\right) \rho
_{\lambda }\left(g', ux_0\right),
\end{equation*}%
or equivalently,%
\begin{equation*}
\rho _{\lambda }\left(g, g'ux_0\right) =\rho _{\lambda }\left(gg', ux_0\right) \rho
_{\lambda }\left(g', ux_0\right) ^{-1}.
\end{equation*}%

The first factor in the right-hand side of the previous equality is estimated by $\rho _{\lambda
}\left(gg', ux_0\right)>c_1$, since $gg'\in \mathrm{int}S$. The second factor satisfies 
$$\rho
_{\lambda }\left(g', ux_0\right) ^{-1}=\rho
_{\lambda }\left(g', x_0\right) ^{-1}\geq M^{-1},\hspace{.5cm}\mbox{ where }\hspace{.5cm}M:=\max_{g'\in U}\rho _{\lambda }\left(g', x_0\right).$$
Since $\pi^{-1}(W)=U\pi^{-1}(\xi_0)=UK_{\Theta}x_0$, we conclude that, 
\begin{equation*}
\rho _{\lambda }\left(g, y\right) >\frac{c_1}{M}>0\qquad y\in \pi^{-1}(W), \;g\in S.
\end{equation*}%

Therefore, by considering $R=\min_{y\in \mathbb{F}}\rho _{\lambda }\left( h,y\right)>0$, we conclude that
\begin{equation*}
\rho _{\lambda }\left(gh, y\right) =\rho _{\lambda }\left( g,hy\right) \rho
_{\lambda }\left( h,y\right)>\frac{c_1R}{M}>0,
\end{equation*}%
proving the result.

\end{prova}

\section{The main result}

In this section we prove our main result, namely we prove the following:

\begin{teorema}
\label{teoprinc}Let $S\subset G$ be a semigroup with $\mathrm{int}S\neq
\emptyset $ and write $\mathbb{F}_{\Theta \left( S\right) }$, $\Theta \left(
S\right) \subset \Sigma $, for its flag type. Denote by $C$ the invariant
control set of $S$ in the maximal flag manifold $\mathbb{F}$.

If $x_{0}\in C_{0}$, then
\begin{equation*}
\inf_{g\in S}\rho _{\alpha }\left( g,x_{0}\right) >0
\end{equation*}%
for any $\alpha \in \Sigma \setminus \Theta \left( S\right) $. Conversely, $%
\inf_{g\in S}\rho _{\alpha }\left( g,x_{0}\right) =0$ if $\alpha \in \Theta
\left( S\right) $.
\end{teorema}

The proof of our main result is divided into the next two sections.

\subsection{Nonexistence of lower bound for $\protect\alpha \in \Theta \left(
S\right) $}

We start by showing that the cocycle $\rho_{\alpha}$, for roots $\alpha$ inside the type flag of $S$, does not admit a lower bound, that is,
\begin{equation*}
\inf_{g\in S}\rho _{\alpha }\left( g,x_{0}\right) =0
\end{equation*}%
if $\alpha \in \Theta \left( S\right) $ and $x_{0}\in C_{0}$.

This is an
immediate consequence of Proposition \ref{propConesEmA}. In fact, let us assume w.l.o.g. that $x_{0}\in C_{0}$ is the origin of $\mathbb{F}$ and
assume further that $A^{+}\cap \mathrm{int}S\neq \emptyset $. By Proposition %
\ref{propConesEmA} it holds that $\Gamma_N\cap r_{\alpha}\mathfrak{a}^+\neq\emptyset$ and hence, there exists $H\in\mathfrak{a}^+$, $t>0$, and $n\in N$ such that $g=\mathrm{e}^{tH}n\in\mathrm{int}S$.

Since $\alpha(H)<0$, we conclude that $\rho
_{\alpha }\left( g,x_{0}\right) =\mathrm{e}^{t\alpha \left( H\right) }<1$, implying that
$$\rho _{\alpha }\left( g^{k},x_{0}\right) =\rho _{\alpha }\left(
g,x_{0}\right) ^{k}\rightarrow 0 \hspace{.5cm}\mbox{ as }\hspace{.5cm} k\rightarrow +\infty.
$$
Therefore, $\rho
_{\alpha }\left( g,x_{0}\right)$ admits no positive lower bound in $S$ when $\alpha\in\Theta(S)$.

\subsection{Simple root outside $\Theta \left( S\right) $}

We apply now the results from Section 3 to get a lower bound for $\rho _{\alpha
}\left(g, x \right) $, $g\in S$, $x\in C_{0}$ when $\alpha \in \Sigma
\setminus \Theta \left( S\right)$. The main idea is to reduce our problem to a rank one group associated with $\alpha$ and apply Corollary \ref{corLowerRankone}.  For this purpose we exploit the
fibration $\mathbb{F}\rightarrow \mathbb{F}_{\alpha}$ over the flag
manifold defined by $\alpha $.

Let $\mathfrak{g}\left( \alpha \right) $ be the real rank one subalgebra
generated by $\mathfrak{g}_{\pm \alpha }$ and write $G\left( \alpha \right) $
for the connected subgroup having Lie algebra $\mathfrak{g}\left( \alpha
\right) $. This subalgebra is the direct sum of $\mathfrak{a}(\alpha )=%
\mathrm{span}\{H_{\alpha }\}$ and the root spaces $\mathfrak{g}_{\pm \alpha
} $ and $\mathfrak{g}_{\pm 2\alpha }$ (eventually $2\alpha $ is not a root
in which case $\mathfrak{g}_{\pm 2\alpha }=\{0\}$). Consider the Iwasawa
decomposition $\mathfrak{g}\left( \alpha \right) =\mathfrak{k}(\alpha)\oplus \mathfrak{a}(\alpha)\oplus \mathfrak{n}(\alpha)$ with $\mathfrak{%
k}(\alpha)=\mathfrak{g}\left( \alpha \right) \cap \mathfrak{k}$ and $\mathfrak{n}(\alpha)=\mathfrak{g}_{\alpha }\oplus \mathfrak{g}_{2\alpha}$
. It defines the global decomposition $G\left( \alpha \right) =K(\alpha)A(\alpha)N(\alpha)$ by taking exponentials of the Lie algebra
components. The restriction of $\alpha$ to $\mathfrak{a}(\alpha)$ is a
root of $\mathfrak{g}\left(\alpha \right) $ which is positive for the
choice of positive roots in $\mathfrak{a}(\alpha)^{\ast }$ that yields
this Iwasawa decomposition.

Take the flag manifold $\mathbb{F}_{\alpha}=G/P_{\alpha }$ and let $%
\pi :\mathbb{F}\rightarrow \mathbb{F}_{\alpha}$ be the canonical
projection. If $x_{0}$ is the origin of $\mathbb{F}=G/P$, then the origin of $%
\mathbb{F}_{\alpha}$ is $x_{\alpha }=\pi \left( x_{0}\right) $. The
fiber $F(\alpha)=\pi ^{-1}\{x_{\alpha }\}$ is invariant by $G\left( \alpha
\right) $ and the action of this group on $F(\alpha)$ is transitive.

The isotropy subalgebra at $x_{0}$ for the action of $G\left( \alpha \right)
$ in $F(\alpha)$ is the parabolic subalgebra $\mathfrak{q}(\alpha)$ of $%
\mathfrak{g}\left( \alpha \right) $ given by $\mathfrak{p}(\alpha)=%
\mathfrak{m}(\alpha)\oplus \mathfrak{a}(\alpha)\oplus \mathfrak{n}(\alpha)$. Hence $F(\alpha)$ becomes identified with the flag manifold
of $G\left( \alpha \right)$, namely, $F(\alpha)=G\left( \alpha
\right) /P(\alpha)$. For this identification, $x_{0}$ is the origin of $F(\alpha)$. Moreover, since $G(\alpha)$ is a real rank one semisimple Lie group, the flag manifold $F(\alpha)$ is diffeomorphic to a sphere (see \cite[Chapter VII, Section 7]{knapp}).
Viewing $\alpha $ as a root in $\mathfrak{a}(\alpha)^{\ast }$ it defines a
cocycle on $G\left( \alpha \right) \times F(\alpha)$. This cocycle is
denoted the same way by $\rho _{\alpha }\left(g, x\right) $ since it is the
restriction to the fiber $F(\alpha)$ of the cocycle over $G\times \mathbb{F%
}$ as follows by the inclusions $K(\alpha)\subset K$, $A(\alpha)\subset
A $ and $N(\alpha)\subset N$.

Therefore,  if $T\subset G\left( \alpha \right) $ is a semigroup with
nonempty interior such that $x_{0}\in C\left( T\right) _{0}$, where $C\left(
T\right) $ is the invariant control set of $T$ in $F(\alpha)$, then $\rho
_{\alpha }\left( g,x_{0}\right) $ is bounded below in $T$ by Corollary \ref%
{corLowerRankone}. This fact will be applied soon to get a lower bound for $%
\rho _{\alpha }\left( g,x_{0}\right) $ in the semigroup $S\subset G$.

As another ingredient for the proof of the required estimate, we consider
decomposition (\ref{eq_langlands}) determined by $\alpha$ given by
\begin{equation*}
P_{\alpha}=MG\left( \alpha \right) A_{\alpha}N_{\alpha}
\end{equation*}%
where $A_{\alpha}=\exp \mathfrak{a}_{\alpha}$, $
N_{\alpha}=\exp \mathfrak{n}_{\alpha}$ and $
\mathfrak{n}_{\alpha} =\sum_{\beta }\mathfrak{g}_{\beta }$ with the sum extended to
the positive roots $\beta \neq \alpha ,2\alpha $.
%is the nilradical of $\mathfrak{p}_{\alpha}$ where

Going back to the semigroup $S$ suppose without loss of generality that $%
x_{0}$ belongs to the core $C_{0}$ of the invariant control set of $S$ in $%
\mathbb{F}$. A semigroup $T\subset G\left( \alpha \right) $ is defined from $%
S$ by the following steps:

\begin{enumerate}
\item The projection $C_{\alpha }=\pi \left( C\right) $ is the invariant
control set of $S$ in $\mathbb{F}_{\alpha}$ whose core $\left(
C_{\alpha }\right) _{0}$ contains $\pi \left( C_{0}\right) $.

\item Let $x_{\alpha }$ be the origin of $\mathbb{F}_{\alpha }$. Then $%
x_{\alpha }=\pi \left( x_{0}\right) \in \pi \left( C_{0}\right) \subset
\left( C_{\alpha }\right) _{0}$. Hence the semigroup $S_{\alpha}=S \cap P_{\alpha}$ has a nonempty interior in $%
P_{\alpha}$.

\item Taking into account the decomposition $P_{\alpha }=MG\left( \alpha
\right) A_{\alpha }N_{\alpha} $ define
\begin{equation*}
\Gamma =\{g\in MG\left( \alpha \right) :\exists b\in A_{\alpha }N_{\alpha } ,~gb\in S\}.
\end{equation*}
Since $MG\left(\alpha \right) $ normalizes $A_{\alpha }N_{\alpha }
$ it follows that $\Gamma $ is a subsemigroup. Moreover, the fact that $S_{\alpha}$ has a nonempty interior in $P_{\alpha}$ implies that $\Gamma$ has a nonempty interior in $%
MG\left( \alpha \right) $.

\item Now define also $T=\Gamma \cap G\left( \alpha \right) $ which is a
subsemigroup of $G\left( \alpha \right) $ with nonempty interior. In
fact, take $mg\in \mathrm{int}\Gamma $, with $m\in M$ and $g\in G\left(
\alpha \right) $ such that $m$ has finite order, say $m^{k}=1$. This is
possible because $M$ is compact so that the set of its elements elements of
finite order is dense. Then $\left( mg\right) ^{k}=m^{k}g_{1}=g_{1}\in
G\left( \alpha \right) $ because $M$ normalizes $G\left( \alpha \right) $.
Thus $g_{1}\in \mathrm{int}\left( \Gamma \cap G\left( \alpha \right) \right)
$, that is, $\mathrm{int}T\neq \emptyset $. 
\end{enumerate}
%The groups $MG\left( \alpha \right) $ and $G\left( \alpha \right) $ act
%transitively on the fiber $F_{\alpha }=\pi ^{-1}\{x_{\alpha }\}$. 
The next
lemma tells about the invariant control sets of $\Gamma $ and $T$ in the
fiber $F(\alpha)$.

\begin{lema}
\label{lemIcsFibra}
Let $C$ be the $S$-invariant control set in $\mathbb{F}$. Then, $F(\alpha)\cap C$ is the unique invariant control set
for $\Gamma $. Moreover $x_{0}$ belongs to the core $C\left( T\right) _{0}$
of the unique invariant control $C\left( T\right) $ of $T$ in $F(\alpha)$.
\end{lema}

\begin{prova}
By general facts of semigroup actions on fiber bundles, $C$ is the union of
invariant control sets on the fibers of $\pi :\mathbb{F}\rightarrow \mathbb{F%
}_{\alpha}$ (see \cite[Theorem 4.4]{Braga}). In particular, $F(\alpha)\cap C$ is the invariant
control set of the subsemigroup $S_{\alpha}$ leaving invariant
the fiber $F(\alpha)$. Now in the decomposition $P_{\alpha}=MG\left(
\alpha \right) A_{\alpha}N_{\alpha} $ if $n\in A_{\alpha}N_{\alpha}$ then its restriction to $F(\alpha)$ is the identity map, since $A_{\alpha}N_{\alpha} $ is a normal
subgroup of $P_{\alpha}$ and $A_{\alpha}N_{\alpha}\subset AN$ implies $nx_0=x_0$. Therefore, the orbits of $\Gamma $ in $F(\alpha)$
are equal to the orbits of $S_{\alpha}$, implying that $F(\alpha)\cap C$ is the invariant control set of $\Gamma $ as well.
Concerning the invariant control set $C\left( T\right) $ of $T$ in $%
F(\alpha)$ we have $C\left( T\right) \subset F(\alpha)\cap C$ because $%
T\subset \Gamma $. Since $x_{0}\in C_{0}$, there exists $g\in \left( \mathrm{%
int}S\right) \cap P$ such that $x_{0}$ is an attractor fixed point of $g$ so
that $g^{k}z\rightarrow x_{0}$, as $k\rightarrow \infty $, for all $z\in
F(\alpha)\cap C$. In view of the decomposition $P_{\alpha }=MG\left(
\alpha \right) A_{\alpha }N_{\alpha }$ and taking into account
that $P\subset P_{\alpha }$, write $g=mhn$ with $m\in M$, $h\in G\left(
\alpha \right) $ and $n\in A_{\alpha }N_{\alpha }$. We have that $%
g_{1}=mh$ belongs to $\mathrm{int}\Gamma $ and $x_{0}$ is an attractor fixed
point of $g_{1}$ in $F(\alpha)$ because $n$ is the identity in this fiber.
We can take $m$ with finite order to conclude that $g_{2}=g_{1}^{k_{0}}\in 
\mathrm{int}T$ for some integer $k_{0}$.

Finally, for any $z\in C\left( T\right) \subset C$, we have $%
g_{2}^{k}z\rightarrow x_{0}$ so that $x_{0}\in C\left( T\right) $. Actually $%
x_{0}$ belongs to the core $C\left( T\right) _{0}$ of $C\left( T\right) $
because $g_{2}\in \mathrm{int}T$ and $g_{2}x_{0}=x_{0}$, concluding the proof.
\end{prova}

Having these constructions at hand, we can get the estimate that allows us to
apply Lemma \ref{lemcompa}.

\begin{lema}
\label{lemcociclolimitadoRaiz}Let $\alpha $ be a simple root outside $\Theta
\left( S\right) $. Then for any $x\in C_{0}$ there exists $d>0$ such that $%
\rho _{\alpha }\left( g,x\right) >d$ for every $g\in S_{x}$.
\end{lema}

\begin{prova}
Assume without loss of generality that $x$ is the origin $x_{0}$ of $\mathbb{%
F}$ so that we can apply the above comments, including Lemma \ref{lemIcsFibra}%
.

Let $d>0$ be a lower bound for $\rho _{\alpha }\left( g,x_{0}\right) $ with $%
g\in T$ which exists by Corollary \ref{corLowerRankone}.

We are required to get such a lower bound for $g$ in $S^{\circ}_{x_{0}}=\left(
\mathrm{int}S\right) \cap P$. We can write $g\in S^{\circ}_{x_{0}}$ as
\begin{equation*}
g=m\left( hn\right) h_{1}n_{1}
\end{equation*}%
with $m\in M$, $h\in A(\alpha)$, $n\in N\left(\alpha\right),$ $%
h_{1}\in A_{\alpha}$ and $n_{1}\in N_{\alpha}$. We have
\begin{equation*}
\rho _{\alpha }\left( g,x_{0}\right) =\rho _{\alpha }\left(
hn,h_{1}n_{1}x_{0}\right) \rho _{\alpha }\left( h_{1}n_{1},x_{0}\right)
=\rho _{\alpha }\left( h,x_{0}\right)
\end{equation*}%
because $m\in M$, $n,n_{1}\in N$, and $\alpha \left( \log h_{1}\right) =0$.
The value of the cocycle does not depend on $m\in M$ so that it can be
changed slightly, keeping $g$ within $\mathrm{int}S$, and assume that $m^{k}=1$
for some integer $k$. In this case if $g_{1}=m\left( hn\right) $ then $%
g_{1}^{k}=h^{k}n^{\prime }$ with $n^{\prime }\in N\left( \alpha \right)$ because $m$ normalizes this subgroup and commutes with $h$. Thus $%
g_{1}^{k}\in \mathrm{int}T$ so that
\begin{equation*}
d<\rho _{\alpha }\left( g_{1}^{k},x_{0}\right) =\rho _{\alpha }\left(
h^{k},x_{0}\right) =\rho _{\alpha }\left( h,x_{0}\right) ^{k}=\rho _{\alpha }\left( g,x_{0}\right) ^{k}.
\end{equation*}%

Therefore we have proved that for any $g\in S_{x_0}^{\circ}=\left( \mathrm{int}S\right) \cap
P$ there exists an integer $k$ such that
\begin{equation*}
\rho _{\alpha }\left( g,x_{0}\right) ^{k}>d.
\end{equation*}%
Hence if $\rho _{\alpha }\left( g,x_{0}\right) <1$ then
\begin{equation*}
\rho _{\alpha }\left( g,x_{0}\right) >\rho _{\alpha }\left( g,x_{0}\right)
^{k}>d
\end{equation*}%
showing that $\min \{d,1\}>0$ is a lower bound for $\rho _{\alpha }\left(
g,x_{0}\right) $ with $g$ running through $\left( \mathrm{int}S\right) \cap
P =S_{x_0}^{\circ}$.
\end{prova}

Combining the above lemma with Lemma \ref{lemcompa}, we get at once the
following statement, which concludes the proof of Theorem \ref{teoprinc}.

\begin{proposicao}
\label{propEstimaRaizFora}For any $\alpha \in \Sigma \setminus \Theta \left(
S\right) $ and $x_0\in C_0$, there exists $c>0$ such that for every $g\in S$ it holds that
\begin{equation*}
\rho _{\alpha }\left( g,x\right) >c.
\end{equation*}
\end{proposicao}

\appendix

\section{Preliminaries}

Here we introduce the basic concepts and results related to semisimple theory, dynamics of semigroups and cocycles.

\subsection{Semisimple theory}

We introduce here the main results and notations related to semisimple theory. For more on the subject, the reader can consult \cite{he1, he2, knapp, knomizuI, warner}.
Let $G$ be a connected, semisimple, non-compact Lie group $G$ with a finite center and associated Lie algebra $\mathfrak{g}$. Choose a Cartan involution $\theta:\mathfrak{g}\rightarrow\mathfrak{g}$ induces an inner product $B_{\theta}(X,Y) = -C(X,\theta(Y))$, where $C(X,Y) = \tr(\ad(X)\ad(Y))$ is the Cartan-Killing form. Since $\theta^2 = \mathrm{id}$, we get that $\theta$ is self-adjoint with respect to $B_{\theta}$ and hence, we have the Cartan decomposition $\mathfrak{g}= \mathfrak{k} \oplus \mathfrak{s}$, where $\mathfrak{k}$ is the $1$-eigenspace and $\mathfrak{s}$ the $(-1)$-eigenspace of $\theta$. We also have an associated Cartan decomposition $G = KS$ of the group with $K = \exp\mathfrak{k}$ and $S = \exp\mathfrak{s}$. The derivations $\ad(X)$, $X\in\mathfrak{k}$, are skew-symmetric w.r.t. $B_{\theta}$ and hence the automorphisms $\Ad(k)$, $k\in K$, are isometries. Let us consider a maximal abelian subspace $\mathfrak{a} \subset\mathfrak{s}$, and for each $\alpha\in\mathfrak{a}^*$ put
\begin{equation*}
 \mathfrak{g}_{\alpha} := \left\{X \in \mathfrak{g}\ :\ \mathrm{ad}(H)X = \alpha(H)X,\ \forall H\in\mathfrak{a}\right\}.
\end{equation*}
The derivations $\ad(H)$, $H\in\mathfrak{a}$, commute and are self-adjoint w.r.t. the inner product $B_{\theta}$. Hence, they can be diagonalized simultaneously, and the nontrivial $\mathfrak{g}_{\alpha}$ are the associated eigenspaces. The set

\begin{equation*}
  \Pi := \left\{\alpha\in\mathfrak{a}^*\backslash\{0\} :\ \mathfrak{g}_{\alpha} \neq 0\right\}%
\end{equation*}
is called the \emph{set of roots of $\mathfrak{g}$}. The associated spaces $\mathfrak{g}_{\alpha}$, $\alpha\in\Pi$, are called \emph{root spaces}. The set of \emph{split-regular elements of $\mathfrak{a}$}\footnote{More generally, any element of the form $\mathrm{Ad}(k)H$ for $H\in\mathfrak{a}^+$ and $k\in K$ is called a split-regular element of $\mathfrak{g}$.} is given by%
\begin{equation*}
  \left\{ H\in\mathfrak{a}\ : \ \alpha(H)\neq0,\ \forall \alpha\in\Pi \right\}.%
\end{equation*}
The connected components of this set are the \emph{Weyl chambers}. Choosing an (arbitrary) Weyl chamber $\mathfrak{a}^+$, we define
\begin{equation*}
  \alpha > 0 \qquad \Leftrightarrow \qquad \alpha|_{\mathfrak{a}^+} > 0,\quad \alpha\in\Pi,%
\end{equation*}
and hence, the \emph{sets of positive and negative roots}, respectively, are
\begin{equation*}
  \Pi^+ := \left\{\alpha\in\Pi\ :\ \alpha>0\right\} \mbox{\quad and \quad} \Pi^- := -\Pi^+,%
\end{equation*}
and $\Pi = \Pi^+ \cup \Pi^-$ is a disjoint union. We define the nilpotent subalgebras $\mathfrak{n} := \sum_{\alpha\in\Pi^+}\mathfrak{g}_{\alpha}$ and $\mathfrak{n}^- := \sum_{\alpha\in\Pi^-}\mathfrak{g}_{\alpha}$, which yields the \emph{Iwasawa decomposition} of the Lie algebra:
\begin{equation*}
  \mathfrak{g} = \mathfrak{k} \oplus \mathfrak{a} \oplus \mathfrak{n}
\end{equation*}
Associated to each of the previous subalgebras are connected Lie subgroups of $G$, denoted by $K$, $A$, $N$, and $N^-$, respectively. They induce the Iwasawa decomposition $G = KAN$ of the group $G$. The \emph{Weyl group} $\mathcal{W}$ of $G$ is the quotient $M^*/M$, where $M^*$ and $M$ are the normalizer and the centralizer of $\mathfrak{a}$ in $K$, respectively, i.e., $M^* = \{k\in K : \Ad(k)\mathfrak{a}=\mathfrak{a}\}$ and $M = \{k\in K : \Ad(k)H=H,\ \forall H\in\mathfrak{a}\}$. 
Alternatively, the Weyl group 
is the group generated by the orthogonal reflections at the hyperplanes $\ker\alpha$, $\alpha\in\Pi$. The Weyl group acts simply transitively on the Weyl chambers. 
%It also acts on the roots by%
%\begin{equation*}
 % w\alpha(H) = \alpha(w^{-1}H),\quad \forall H\in\mathfrak{a},\ \alpha\in\Pi.%
%\end{equation*}
The \emph{set of simple roots}  $\Sigma \subset \Pi^+$ are the positive roots that cannot be written as linear combinations of other positive roots. It forms a basis of $\mathfrak{a}^*$. Moreover, the reflections at $\ker\alpha$, $\alpha\in\Sigma$, generate the Weyl group. There exists a unique element $w_0\in\mathcal{W}$ of order $2$ which takes $\Pi^+$ to $\Pi^-$, called the \emph{principal involution}. 
Let $\Theta \subset \Sigma$ be an arbitrary subset. We write $\langle\Theta\rangle$ for the set of roots that are linear combinations (over $\mathbb{Z}$) of elements in $\Theta$. Moreover, we put $\mathfrak{a}(\Theta) := \langle H_{\alpha}\ : \ \alpha\in\Theta \rangle$, where $H_{\alpha}\in\mathfrak{a}$ is the coroot of $\alpha$, defined by $B_{\theta}(H_{\alpha},H) = \alpha(H)$, $H\in\mathfrak{a}$. The subalgebra $\mathfrak{g}(\Theta)$ generated by $\mathfrak{a}(\Theta) \oplus \sum_{\alpha\in\langle\Theta\rangle}\mathfrak{g}_{\alpha}$ is a semisimple Lie algebra. We put $\mathfrak{k}(\Theta) := \mathfrak{k} \cap \mathfrak{g}(\Theta)$ and $\mathfrak{n}(\Theta) := \mathfrak{n} \cap \mathfrak{g}(\Theta)$. Then $\mathfrak{g}(\Theta)$ is the Lie algebra of a semisimple Lie group $G(\Theta)\subset G$, and $\mathfrak{g}(\Theta) = \mathfrak{k}(\Theta)\oplus\mathfrak{a}(\Theta)\oplus\mathfrak{n}(\Theta)$ is an Iwasawa decomposition of $\mathfrak{g}(\Theta)$, while $\Theta$ is the corresponding set of simple roots. We write $K(\Theta)$ for the connected Lie subgroup with Lie algebra $\mathfrak{k}(\Theta)$ and $A(\Theta) = \exp\mathfrak{a}(\Theta)$, $N(\Theta) = \exp\mathfrak{n}(\Theta)$ (which are also connected subgroups). Then $G(\Theta) = K(\Theta)A(\Theta)N(\Theta)$ is an Iwasawa decomposition of $G(\Theta)$. Let
\begin{equation*}
  \mathfrak{a}_{\Theta} := \left\{ H\in\mathfrak{a}\ :\ \alpha(H)=0,\ \forall \alpha \in \Theta \right\}%
\end{equation*}
be the orthogonal complement of $\mathfrak{a}(\Theta)$ in $\mathfrak{a}$.
%and{\color{blue}  note that%
%\begin{equation*}
 % \mathfrak{a}(\Theta_1 \cap \Theta_2) = \mathfrak{a}(\Theta_1) \cap \mathfrak{a}(\Theta_2) \mbox{\quad and\quad } \mathfrak{a}_{\Theta_1 \cap \Theta_2} = \mathfrak{a}_{\Theta_1} + \mathfrak{a}_{\Theta_2}.
%\end{equation*}
%}
The subset $\Theta$ singles out a subgroup $\mathcal{W}_{\Theta}$ of $\mathcal{W}$ consisting of those elements that act trivially on $\mathfrak{a}_{\Theta}$. Alternatively, $\mathcal{W}_{\Theta}$ can be defined as the subgroup generated by the reflections at the hyperplanes $\ker\alpha$, $\alpha\in\Theta$. Then $\mathcal{W}_{\Theta}$ is isomorphic to the Weyl group $\mathcal{W}(\Theta)$ of $G(\Theta)$. We let $Z_{\Theta}$ denote the centralizer of $\mathfrak{a}_{\Theta}$ in $G$ and $K_{\Theta} = Z_{\Theta} \cap K$
implying
%\begin{equation}\label{eq_kthetaintersec}
% K_{\Theta_1 \cap \Theta_2} = K_{\Theta_1} \cap K_{\Theta_2}.%
%\end{equation}
The group $Z_{\Theta}$ is a reductible Lie group and admits an Iwasawa decomposition that is given by
$Z_{\Theta} = K_{\Theta}AN(\Theta).$
The \emph{parabolic subalgebra of type $\Theta\subset\Sigma$} is defined by $\mathfrak{p}_{\Theta} := \mathfrak{n}^-(\Theta) \oplus \mathfrak{p}$, where $\mathfrak{p} := \mathfrak{m}\oplus \mathfrak{a}\oplus \mathfrak{n}$, $\mathfrak{n}^-(\Theta)=\mathfrak{n}^-\cap\mathfrak{g}(\Theta)$  , and $\mathfrak{m}$ is the Lie algebra of $M$, the centralizer of $\mathfrak{a}$ in $\mathfrak{k}$, or, respectively, the part of the common $0$-eigenspace of the maps $\ad(H)$, $H\in\mathfrak{a}$, contained in $\mathfrak{k}$. The associated \emph{parabolic subgroup} $P_{\Theta}$ is the normalizer of $\mathfrak{p}_{\Theta}$ in $G$. Then $\mathfrak{p}_{\Theta}$ is the Lie algebra of $P_{\Theta}$. The empty set $\Theta=\emptyset$ yields the minimal parabolic subalgebra $\mathfrak{p}_{\emptyset} = \mathfrak{p}$. An Iwasawa decomposition of the associated subgroup $P$ is $P = MAN^+$. The parabolic subgroup $P_{\Theta}$ also decomposes as
\begin{equation}\label{eq_langlands}
  P_{\Theta} = K_{\Theta}AN=MG(\Theta)A_{\Theta}N_{\Theta},
\end{equation}
where $A_{\Theta}=\exp\mathfrak{a}_{\Theta}$ and $N_{\Theta}=\exp\,\mathfrak{n}_{\Theta}$, with $\mathfrak{n}_{\Theta}=\sum_{\alpha\in\Pi^+\setminus\langle\Theta\rangle}\mathfrak{g}_{\alpha}$ the nilradical of $\mathfrak{p}_{\Theta}$. The first equality of (\ref{eq_langlands}) is known as \emph{Langland decomposition}, and the second one follows from the fact that $K_{\Theta}=MK(\Theta)$, $A=A_{\Theta}A(\Theta)=A(\Theta)A_{\Theta}$, and $N=N(\Theta)N_{\Theta}=N_{\Theta}N(\Theta)$.

The \emph{flag manifold of type $\Theta$} is the $\Ad(G)$-orbit $\mathbb{F}_{\Theta} := \Ad(G)\mathfrak{p}_{\Theta}$ with base point $x_{\Theta} := \mathfrak{p}_{\Theta}$ in the Grassmann manifold of $(\dim\mathfrak{p}_{\Theta})$-dimensional subspaces of $\mathfrak{g}$. $\mathbb{F}_{\Theta}$ is compact, because the natural $K$-action on $\mathbb{F}_{\Theta}$ is transitive; in fact, it holds that $\mathbb{F}_{\Theta}=K/K_{\Theta}$. In the case $\Theta=\emptyset$ we also write $x_0 =x_{\emptyset}$ and $\mathbb{F} = \mathbb{F}_{\emptyset}$ for the maximal flag manifold. Since the isotropy group of $x_{\Theta}$ is the subgroup $P_{\Theta}$, $\mathbb{F}_{\Theta}$ can be identified with the homogeneous space $G/P_{\Theta}$ via $\Ad(g)\mathfrak{p}_{\Theta}\mapsto gP_{\Theta}$. If $\Theta_1 \subset \Theta_2$, then $P_{\Theta_1}\subset P_{\Theta_2}$ and the projection $\pi^{\Theta_1}_{\Theta_2}:\mathbb{F}_{\Theta_1} \rightarrow \mathbb
{F}_{\Theta_2}$, $gP_{\Theta_1} \mapsto gP_{\Theta_2}$, is a well-defined fibration. In case $\Theta_1=\emptyset$, we just write $\pi_{\Theta_2}$ for this map.
The choice of the previous subalgebras and subgroups is not unique. In fact, by conjugation with any $k\in K$, one obtains a new set of subalgebras and subgroups. For instance, $\mathrm{Ad}(k)\mathfrak{a}$ is another maximal abelian subspace of $\mathfrak{s}$, and $\mathrm{Ad}(k)\mathfrak{a}^+$ is an associated positive Weyl chamber. Also $\mathfrak{g}= \mathfrak{k} \oplus \mathrm{Ad}(k)\mathfrak{a} \oplus \mathrm{Ad}(k)\mathfrak{n}$ is another Iwasawa decomposition. We will use these conjugated settings frequently.%

\subsection{Semigroup actions on flag manifolds}

In this section we highlight the main results involving the actions of semigroups on the flag manifolds. The results presented here can be found in the references \cite{smmax, smc, smorder, smt}.

Let $S\subset G$ be a semigroup with nonempty interior. The semigroup $S$ acts on
the flag manifolds $\mathbb{F}_{\Theta}$ of $G$ by left translation. This action is $S$ is
not transitive in $\mathbb{F}_{\Theta}$ unless $S=G$.

A \emph{control set} for the $S$-action on $\mathbb{F}_{\Theta}$ is a subset $D\subset\mathbb{F}_{\Theta}$ that is maximal w.r.t. set inclusion, satisfying:
\begin{enumerate}
 \item $\mathrm{int}D\neq\emptyset$;
\item $D\subset\mathrm{cl}(Sx)$ for all $x\in D$.
\end{enumerate}
We say that $D$ is $S$-invariant, provided that $SD\subset D$.

As shown in \cite[Theorem 3.1]{smc}, there exists exactly one invariant control set $C_{\Theta}\subset\mathbb{F}_{\Theta}$, with $C_{\Theta}\neq\mathbb{F}_{\Theta}$ if $S$ is a proper semigroup. The \emph{core} of $C_{\Theta}$ is the open and dense subset given by
$$(C_{\Theta})_0:=\{x\in C_{\Theta}; x\in\mathrm{int}Sx\}.$$
Let $H\in\mathfrak{a}^+$ and $k\in K$, and consider a split-regular element $Z=\Ad(k)H$. The action of $Z$ on the flag manifold $\mathbb{F}_{\Theta}$ admits as fixed points the elements of the form $\mathrm{fix}_{\Theta}(Z, w)=kwx_{\Theta}$ for $w\in\mathcal{W}$. Those points are isolated and have stable and unstable manifolds given by $\mathrm{st}_{\Theta}(Z, w):=kN^-wx_{\Theta}$ and $\mathrm{un}_{\Theta}(Z, w)=kNwx_{\Theta}$, respectively. In particular, there is a unique attractor fixed point $\mathrm{at}_{\Theta}(Z)=kx_{\Theta}$, whose stable manifold $\mathrm{st}_{\Theta}(Z)=\mathrm{st}_{\Theta}(Z, 1)=kN^-x_{\Theta}$ is open and dense and a unique repeller fixed point $\mathrm{rp}_{\Theta}=kw_0x_{\Theta}$ whose unstable manifold $\mathrm{un}_{\Theta}(Z)=\mathrm{st}_{\Theta}(Z, w_0)=kN^+w_0x_{\Theta}$ is also open and dense. In particular, the core of the invariant control set $C_{\Theta}$ is characterized by attractor points \cite[Theorem 3.4]{smt} as
$$(C_{\Theta})_0=\{\mathrm{at}_{\Theta}(Z), \;Z \mbox{ is split-regular and }\mathrm{e}^Z\in\mathrm{int} S \}$$

The geometry of the invariant control sets associated with a semigroups in $G$ is described  by the following results, proved in \cite{smt}.

\begin{proposicao}
There exists $\Theta\subset\Sigma$ such that $\pi^{-1}_{\Theta}(C_{\Theta})=C$ is the invariant control set in the maximal flag manifold $\mathbb{F}$. Among the subsets $\Theta$ satisfying
this property there exists a unique maximal one (w.r.t. set inclusion). 
\end{proposicao}

The maximal subset of the previous result is called the \emph{flag type} of $S$ and is denoted by $\Theta(S)$.  Alternatively, we call the flag type of $S$ the corresponding flag manifold $\mathbb{F}_{\Theta(S)}$ (see \cite{smmax, smorder, smt, smmoment} for further
discussions about the flag type of a semigroup and its applications).

Another characterization of the flag type is given by the following result

\begin{proposicao}
    The flag type $\Theta(S)$ is the minimal subset (w.r.t. set inclusion) satisfying: If $Z\in\mathfrak{g}$ is split-regular, then $$\mathrm{e}^Z\in\mathrm{int}S\hspace{.5cm}\implies\hspace{.5cm}C_{\Theta}\subset\mathrm{st}_{\Theta}(Z).$$
\end{proposicao}

\subsection{Cocycles on flag manifolds}

A cocycle $\rho$ on the maximal flag manifold $\mathbb{F}$ is a function $\rho:G\times\mathbb{F}\rightarrow\mathbb{R},$
satisfying
$$\rho(gh, x)=\rho(g, hx)\rho(h, x)\hspace{.5cm}\forall g, h\in G, x\in\mathbb{F}.$$
A natural cocycle on $\mathbb{F}$  appears naturally from an Iwasawa decomposition. In fact, let us start with the map $\mathsf{a}:G\times K\rightarrow \mathfrak{a}$ define by the Iwasawa decomposition
$$gk=u\exp(\mathsf{a}(g, k))n,\hspace{.5cm} u\in K, n\in N.$$
Since the second component of $\mathsf{a}$
is invariant by $M$, the map $\mathsf{a}$ factors to the flag $\mathbb{F}$. Since $\mathsf{a}$ satisfies
$$\mathsf{a}(gh, x)=\mathsf{a}(g, hx)+\mathsf{a}(h, x),$$
it holds that, for any $\lambda\in\mathfrak{a}$, the map $$\rho_{\lambda}:G\times\mathbb{F}\rightarrow\mathbb{R}, \hspace{.5cm}\rho_{\lambda}(g, x)=\mathrm{e}^{\lambda\mathsf{a}(g, x)},$$
is a cocycle\footnote{In fact, any $K$-invariant cocycle on $\mathbb{F}$ is of the previous form for some $\lambda\in\mathfrak{a}^*$.}. The next result \cite[Lemma 6.1]{asm} shows that a cocycle can factors to the partial flag manifolds under some conditions on the functional $\lambda$.

\begin{lema}
\label{restriction}
    If $\lambda\in \mathfrak{a}$ satisfies $\lambda(\mathfrak{a}(\Theta))=0$, then 
    $$\rho_{\lambda
}\left(g, x\right)=\rho_{\lambda
}\left(gu, x\right), \hspace{.5cm}\forall u\in K_{\Theta}.$$
\end{lema}

\end{document}